\documentclass[11pt]{article}
\usepackage{}
\usepackage{bbm}
\usepackage{amssymb}
\usepackage{mathrsfs}%有更多的数学符号，比如花体$\mathscr P$
\usepackage{tikz}
\usepackage{mathrsfs}
\usepackage{amssymb}
\usepackage{amsmath}
\usepackage{bm}
\usepackage{graphicx}
\usepackage{pst-poly}  %画图
\usepackage{pst-plot}  %坐标
\usepackage{url}  %\url{输入你的网址，直接用~即可}

\newtheorem{thm}{Theorem}[section]
\newtheorem{lem}[thm]{Lemma}

\newtheorem{rem}[thm]{Remark}

\newenvironment{pf}[1][Proof]{\noindent\textbf{#1.} }{\hfill\rule{1mm}{2mm}}

\makeatletter \@addtoreset{equation}{section} \makeatother

\begin{document}

\title
{Sufficient conditions for graphs to be $k$-connected, maximally connected and super-connected
\thanks{This work was supported by the National Natural Science Foundation of China (No. 11601002, 11601001)
and the University Natural Science Research Project of Anhui Province (No. KJ2016A003).}}

\author
{Zhen-Mu Hong\thanks{Corresponding author},\ \ Zheng-Jiang Xia \vspace{2mm}\\
{\small School of Finance, Anhui University of Finance and Economics, 233030 Bengbu, China}\\
{\small {\tt zmhong@mail.ustc.edu.cn, xzj@mail.ustc.edu.cn}}\\
\\
Fuyuan Chen \vspace{2mm}\\
{\small Institute of Statistics and Applied Mathematics, }\\
{\small Anhui University of Finance and Economics, 233030 Bengbu, China}\\
{\small {\tt chenfuyuan19871010@163.com}}\\
\\
Lutz Volkmann \vspace{2mm}\\
{\small Lehrstuhl II f\"{u}r Mathematik, RWTH Aachen University, 52056 Aachen, Germany}\\
{\small {\tt volkm@math2.rwth-aachen.de}}\\}

%\\
%Jun-Ming Xu\\
%{\small School of Mathematical Sciences}\\
%{\small University of Science and Technology of China}\\
%{\small Hefei, Anhui, 230026, China}\\
%{\small Email: xujm@ustc.edu.cn}\\}

\date{}
\maketitle

\begin{center}
\begin{minipage}{140mm}

{\bf Abstract:} Let $G$ be a connected graph with minimum degree $\delta(G)$ and vertex-connectivity $\kappa(G)$.
The graph $G$ is $k$-connected if $\kappa(G)\geq k$,
maximally connected if $\kappa(G) = \delta(G)$, and super-connected (or super-$\kappa$) if every minimum vertex-cut
isolates a vertex of minimum degree. In this paper, we show that a connected graph or a connected triangle-free graph
is $k$-connected, maximally connected or super-connected if the number of edges or the spectral radius is large enough.

\vskip16pt



\end{minipage}
\end{center}

\section{Introduction}

Let $G=(V, E)$ be a simple connected undirected graph, where $V=V(G)$ is the
vertex-set of $G$ and $E=E(G)$ is the edge-set of $G$.
The {\it order} and {\it size} of $G$ are defined by $n = |V(G)|$ and $m = |E(G)|$, respectively;
$d_G(x)$ is the {\it degree} of a vertex $x$ in $G$, i.e. the number of edges incident with
$x$ in $G$; $\delta(G)=\min\{d_G(x): x\in V(G)\}$ is the {\it minimum
degree} of $G$. For a subset $X\subset V (G)$,
use $G[X]$ to denote the subgraph of $G$ induced by $X$.
For two subsets $X$ and $Y$ of $V(G)$, let $[X, Y]$ be the set of edges between $X$ and $Y$.
The {\it complement} of $G$ is denoted by $\overline{G}$.
%We write $K_n$ for the complete graph of order $n$.
Let $G_1\cup G_2$ denote the disjoint union of graphs $G_1$ and $G_2$, and
$G_1\vee G_2$ denote the graph obtained from $G_1\cup G_2$ by joining each vertex of $G_1$ to each vertex of $G_2$.
Denote by $\rho(G)$ the largest eigenvalue or the spectral radius of the adjacency matrix of $G$ and
called it the {\it spectral radius} of $G$. If $G$ is connected, then by Perron-Frobenius Theorem,
$\rho(G)$ is simple and there exists a unique (up to a multiple) corresponding positive eigenvector.

A {\it vertex-cut} of a connected graph $G$ is a set of
vertices whose removal disconnects $G$. The {\it vertex-connectivity} or simply the {\it connectivity} $\kappa=\kappa(G)$ of a connected
graph $G$ is the minimum cardinality of a vertex-cut of $G$ if $G$ is not complete,
and $\kappa(G)=n-1$ if $G$ is the complete graph $K_n$ of order $n$.
A vertex-cut $S$ is a {\it minimum vertex-cut} or a {\it $\kappa$-cut} of $G$ if $|S|=\kappa(G)$. Apparently,
$\kappa(G)\leq \delta(G)$ for any graph $G$.
The graph $G$ is {\it $k$-connected} if $\kappa(G)\geq k$,
{\it maximally connected} if $\kappa(G) = \delta(G)$, and {\it super-connected} (or {\it super-$\kappa$})
if every minimum vertex-cut isolates a vertex of minimum degree.
Hence, every super-connected graph is also maximally connected.
An {\em edge-cut} of a connected graph $G$ is a set of
edges whose removal disconnects $G$. The {\em edge connectivity} $\lambda=\lambda(G)$ of a connected
graph $G$ is defined as the minimum cardinality of an edge-cut over all edge-cuts of $G$. An edge-cut
$S$ is a {\em minimum edge-cut} if $|S|=\lambda(G)$. The inequality
$\lambda(G)\le\delta(G)$ is obvious. The graph $G$ is {\em maximally edge-connected} if
$\lambda(G)=\delta(G)$, and it is {\em super-edge-connected} if every minimum edge-cut
consists of edges incident with a vertex of minimum degree.
Therefore, every super-edge-connected graph is also maximally edge-connected.
For graph-theoretical terminology and notation not defined here,
one can refer to \cite{x01,x03}.

%In this paper, we are interested in the conditions to guarantee graphs to be maximally connected or super-connected.

Sufficient conditions for graphs to be maximally (edge-)connected or super-(edge-)connected
were given by several authors, depending on the order, the maximum and minimum degree, the diameter,
the girth, the degree sequence, the clique number and so on. The paper \cite{hv08} by
Hellwig and Volkmann gives a survey on this work. Recently, Volkmann and Hong \cite{vh17}
proved that a connected graph or a connected triangle-free graph is
maximally edge-connected or super-edge-connected if the number of edges is
large enough. Meanwhile, we notice that the relationship between
graph properties and eigenvalues has attracted much attention,
especially the spectral conditions for graphs to be hamiltonian or traceable
since the publication of the paper by Fiedler and Nikiforov \cite{fn10}.
Their paper lays a foundation for the further study to this field. Soon afterwards,
using spectral radius and signless Laplacian spectral radius,
Li \cite{li14} presented sufficient conditions for a graph to be $k$-connected;
Feng et. al. \cite{fzll17} demonstrated sufficient conditions based on spectral radius
for a graph to be $k$-connected and $k$-edge-connected; Feng et. al. \cite{fzl17}
gave a tight sufficient condition for a connected graph with fixed minimum degree
to be $k$-connected based on its spectral radius, for sufficiently large order.
Motivated by the researches above, this paper will prove that a connected graph or a connected triangle-free graph
is $k$-connected, maximally connected or super-connected if the number of edges or the spectral radius is large enough.
The results in this paper improve the conclusion in the paper \cite{fzl17} by Feng et. al. to some extent.

The rest of this paper is organized as follows. Section 2 presents sufficient conditions
for graphs to be $k$-connected depending on size (i.e. the number of edges) or the spectral radius of graph and its complement.
In Section 3, by setting $k=\delta$, we get sufficient conditions for graphs to be maximally connected
depending on size or the spectral radius of graph and its complement.
In Section 4, we obtain sufficient conditions for graphs to be super-connected
depending on size or the spectral radius of graph and its complement.
In the last Section, sufficient conditions for triangle-free graphs to be
$k$-connected, maximally connected or super-connected are acquired depending on the number of edges.

%However, there seems to be no result in the literature showing that graphs of sufficiently large size
%or large spectral radius are maximally connected. In this paper, we discuss the sufficient conditions
%for graphs to be maximally connected or super-connected from this point of view.
%
%
%We show that a connected graph $G$ is maximally connected or super-$\kappa$ if the number of edges is large enough
%or the spectral radius is large enough. For example, for any connected graph $G$, if
%$m \geq {n-2 \choose 2}+2\delta-1$, or $n\geq \delta^2-2\delta+7$ and $\rho(G) \geq n-3$,
%then $G$ is maximally connected, unless $G=K_{\delta-1}\vee (K_2\cup K_{n-\delta-1})$.

\section{$k$-connected graphs}

%\begin{lem}{\rm (Hong~\cite{h88})}\label{lem1.1}
%Let $G$ be a connected graph with $n$ vertices and $m$ edges. Then
%$$\rho(G)\leq \sqrt{2m-n+1}$$
%and the equality holds if and only if $G\cong K_n$ or $G\cong K_{1,n-1}$.
%\end{lem}

For any connected graph $G$ of order $n$ and minimum degree $\delta$,
if $n\leq 4$ or $\delta=1$, then $\kappa=\delta$.
If $\delta=n-1$, then $G = K_n$ and $\kappa=\delta$. If $\delta=n-2$,
then when $u$ and $v$ are nonadjacent the other $n-2$ vertices are all common
neighbors of $u$ and $v$. It is necessary to delete all common neighbors of
some pair of vertices to separate the graph, so $\kappa\geq n-2=\delta$.
Therefore, we only need to consider $n\geq 5$ and $2\leq\delta\leq n-3$ in the following.

\begin{thm}\label{thm2.1}
Let $G$ be a connected graph of order $n\geq 5$, size $m$, minimum degree $\delta\geq k\geq  2$ and vertex-connectivity $\kappa$.

{\rm (a)}\
If
\begin{equation}\label{e2.1}
m \geq \frac12 n(n-1) - (\delta-k+2)(n-\delta-1),
\end{equation}
then $\kappa\geq k$, unless $G=K_{k-1}\vee (K_{\delta-k+2}\cup K_{n-\delta-1})$.

{\rm (b)}\
If
$n\geq \frac12(k+1)(\delta-k+2)+(\delta+2)$ and
\begin{equation}\label{e2.2}
%m \geq \frac12 n(n-1) - (\delta-k+2)(n-\delta-1)-\frac12(k-1)(\delta-k+2),
m \geq \frac12 n(n-1) - \frac12(\delta-k+2)(2n-2\delta+k-3),
\end{equation}
then $\kappa\geq k$, unless $G$ is a subgraph of $K_{k-1}\vee (K_{\delta-k+2}\cup K_{n-\delta-1})$.
\end{thm}

\begin{pf}
Suppose that $1\leq \kappa\leq k-1$.
%It suffices to prove $G=K_{\delta-1}\vee (K_2\cup K_{n-\delta-1})$.
Let $S$ be an arbitrary minimum vertex-cut, and
let $X_0, X_1,\dots$, $X_{p-1}$ ($p\geq 2$) denote the vertex sets of the components of $G-S$, where
$|X_0|\leq |X_1|\leq \cdots\leq |X_{p-1}|$. Each vertex in
$X_i$ is adjacent to at most $|X_i|-1$ vertices of $X_i$ and $\kappa=|S|$ vertices of $S$. Thus
\begin{equation*}
\delta|X_i| \leq \sum_{x\in X_i}d(x) \leq |X_i|(|X_i|+\kappa-1),
\end{equation*}
and so
$\delta-\kappa+1 \leq |X_i| \leq n-\kappa-(\delta-\kappa+1),$
which means
\begin{equation}\label{e2.3}
\delta-\kappa+1 \leq |X_i| \leq n-\delta-1 \text{~for~} i=0,1,\dots,p-1.
\end{equation}
Let $Y=\bigcup_{i=1}^{p-1} X_i$. Then, by (\ref{e2.3}), $\delta-\kappa+1 \leq |Y| \leq n-\delta-1$.
Since $G-S$ is disconnected, there are no edges between $X_0$ and $Y$ in $G$ and
\begin{equation}\label{e2.4}
m \leq \frac12 n(n-1) - |X_0|\cdot|Y|.
\end{equation}

(a) It suffices to prove that $G=K_{k-1}\vee (K_{\delta-k+2}\cup K_{n-\delta-1})$.
By (\ref{e2.3}) and $|X_0|+|Y|=n-\kappa$, we obtain
\begin{align*}
m &\leq \frac12 n(n-1) - |X_0|\cdot|Y| \\
  &\leq \frac12 n(n-1) - (\delta-\kappa+1)(n-\delta-1) ~~~(\text{as}~\kappa\leq k-1)\\
  &\leq \frac12 n(n-1) - (\delta-k+2)(n-\delta-1).
\end{align*}
Combining this with (\ref{e2.1}), we obtain $m=\frac12 n(n-1) - (\delta-k+2)(n-\delta-1)$.
Hence, all the inequalities in the proof above must be equalities.
Thus, we have $\kappa=k-1$, $p=2$, $|X_0|=\delta-k+2$,
$|Y|=n-\delta-1$, $d_G(s)=n-1$ for each $s\in S$,
$d_{G[X_0]}(x)=\delta-k+1$ for each $x\in X_0$ and $d_{G[Y]}(y)=n-\delta-2$ for each $y\in Y$. That is
$G[X_0]=K_{\delta-k+2}$, $G[S]=K_{k-1}$, $G[Y]=K_{n-\delta-1}$.
It follows that $G=K_{k-1}\vee (K_{\delta-k+2}\cup K_{n-\delta-1})$.

(b)
To prove $G$ is a subgraph of $K_{k-1}\vee (K_{\delta-k+2}\cup K_{n-\delta-1})$, it suffices to show that $|X_0|=\delta-k+2$.
If $|X_0|\geq \delta-k+3$, then, combining $|X_0|+|Y|=n-\kappa$ with (\ref{e2.2}) and (\ref{e2.4}), we obtain
\begin{align*}
\frac12 n(n-1) - \frac12(\delta-k+2)(2n-2\delta+k-3) \leq
m & \leq  \frac12 n(n-1) - |X_0|\cdot|Y| \\
&\leq  \frac12 n(n-1)-(\delta-\kappa+2)(n-\delta-2) \\ %\  \  (\text{as~} \kappa\leq k-1)\\
&\leq  \frac12 n(n-1)-(\delta-k+3)(n-\delta-2),
\end{align*}
which yields $n\leq \frac12(k+1)(\delta-k+2)+(\delta+2)$, and so $n = \frac12(k+1)(\delta-k+2)+(\delta+2)$.
It is easy to verify that $\kappa=k-1$, $|X_0|=\delta-k+3$, $|Y|=n-\delta-2$
and $G=K_{k-1}\vee (K_{\delta-k+3}\cup K_{n-\delta-2})$. However,
$\delta(G)=\delta(K_{k-1}\vee (K_{\delta-k+3}\cup K_{n-\delta-2}))=\delta+1>\delta$, which is a contradiction.
Thus, $|X_0|\leq \delta-k+2$. Combining this with $|X_0|\geq \delta-\kappa+1\geq \delta-k+2$,
we get $|X_0|=\delta-k+2$. Since $|S|=\kappa=k-1$ and $d_G(x)\geq \delta$
for each $x\in X_0$, we have each vertex of
$X_0$ is adjacent to each vertex of $S$. Therefore,
$G$ is a subgraph of $K_{k-1}\vee (K_{\delta-k+2}\cup K_{n-\delta-1})$.
\end{pf}

%\begin{thm}\label{thm2.3}
%Let $G$ be a connected graph of order $n\geq 5$, minimum degree $\delta\geq 2$ and vertex-connectivity $\kappa$. If
%\begin{equation}\label{e2.5}
%\rho(G) \geq \frac{\delta-1}{2}+\sqrt{(n-\delta-1)(n-4) + \frac{(\delta+1)^2}{4}},
%\end{equation}
%then $\kappa=\delta$, unless $G=K_{n-4}\vee (K_2\cup K_2)$.
%\end{thm}
%
%\begin{pf}
%Since $G$ is connected, by (\ref{e2.5}) and Lemma~\ref{lem2.4}, we have
%\begin{equation}\label{e2.6}
%\frac{\delta-1}{2}+\sqrt{(n-\delta-1)(n-4) + \frac{(\delta+1)^2}{4}}  \leq \rho(G) \leq \frac{\delta-1}{2}+\sqrt{2|E(G)|-\delta n+\frac{(\delta+1)^2}{4}},
%\end{equation}
%which yields
%$$|E(G)|\geq {n-2 \choose 2}+2\delta-1.$$
%Suppose that $\kappa<\delta$.
%By Theorem~\ref{thm2.1} (a), $G=K_{\delta-1}\vee (K_2\cup K_{n-\delta-1})$.
%To complete the proof, we only need to show $\delta=n-3$.
%
%Since $|E(G)|={n-2 \choose 2}+2\delta-1$, the equalities hold in (\ref{e2.6}).
%Thus, by Lemma~\ref{lem2.4}, $G$ is a bidegreed graph in which each vertex is of degree either $\delta$ or $n-1$.
%However, the vertices of $G$ have degrees from the set $\{\delta,n-3,n-1\}$.
%Therefore, $\delta=n-3$ and the result follows.
%\end{pf}

\begin{thm}\label{thm2.2}
Let $G$ be a connected graph of order $n$, minimum degree $\delta\geq k\geq 2$ and vertex-connectivity $\kappa$. If
\begin{equation}\label{e2.5}
\rho(G) \geq \rho(K_{k-1}\vee (K_{\delta-k+2}\cup K_{n-\delta-1})) ,
\end{equation}
then $\kappa\geq k$, unless $G=K_{k-1}\vee (K_{\delta-k+2}\cup K_{n-\delta-1})$, where $\rho(K_{k-1}\vee (K_{\delta-k+2}\cup K_{n-\delta-1}))$
is the largest root of the equation
$$\lambda^3 - (n - 3)\lambda^2 + ((\delta-k+2)(n-\delta-1)-2n+3)\lambda +(\delta-k+2)(n-\delta-1)k - n+1=0.$$
\end{thm}

\begin{pf}
Assume that (\ref{e2.5}) holds but $1\leq \kappa\leq k-1$.
Let $S$ be an arbitrary minimum vertex-cut of $G$, and
let $X_0, X_1,\dots$, $X_{p-1}$ ($p\geq 2$) denote the vertex sets of the components of $G-S$, where
$|X_0|\leq |X_1|\leq \cdots\leq |X_{p-1}|$. Each vertex in
$X_i$ is adjacent to at most $|X_i|-1$ vertices of $X_i$ and $\kappa=|S|$ vertices of $S$. Thus
\begin{equation*}
\delta|X_i| \leq \sum_{x\in X_i}d(x) \leq |X_i|(|X_i|-1+\kappa) ,
\end{equation*}
and so
$|X_i|\geq \delta-\kappa+1 $ for each $i=0,1,\dots,p-1$.
Let $Y = \bigcup_{i=1}^{p-1} X_i$. Then $\delta-\kappa+1 \leq |X_0|\leq |Y| \leq n-\delta-1$ and $|X_0|+|Y| = n-\kappa$.
Since there are no edges between $X_0$ and $Y$ in $G$,
$G$ is a subgraph of $K_{\kappa}\vee (K_{|X_0|}\cup K_{|Y|})$ and
$\rho(G)\leq \rho(K_{\kappa}\vee (K_{|X_0|}\cup K_{|Y|}))$.

Next, we will show
$$\rho(K_{\kappa}\vee (K_{|X_0|}\cup K_{|Y|}))\leq \rho(K_{\kappa}\vee (K_{\delta-\kappa+1}\cup K_{n-\delta-1}))\leq \rho(K_{k-1}\vee (K_{\delta-k+2}\cup K_{n-\delta-1})).$$

In short, denote $G(a,b,\kappa)=K_{\kappa}\vee (K_a\cup K_b)$ where $b\geq a\geq \delta-\kappa+1$ and $a+b+\kappa=n$.
Let $\textbf{x} = (x_1, x_2,\dots , x_n)^T$ be the unique positive unit eigenvector
corresponding to $\rho(G(a,b,\kappa))$. By symmetry,
let $x:= x_i $ for any $i \in K_a$; $y:= x_j$ for any
$j\in K_\kappa$; $z := x_\ell $ for any $\ell \in K_b$.
According to $\lambda x_i=\sum\limits_{ij\in E(G(a,b,\kappa))}x_j$ and the uniqueness of $\textbf{x}$, we have that
$\rho(G(a,b,\kappa))$ is the largest root of following equations:
\begin{align*}
\lambda x &= (a-1)x+\kappa y,\\
\lambda y &= ax+(\kappa-1)y+bz,\\
\lambda z &= \kappa y+(b-1)z.
\end{align*}
Thus $\rho(G(a,b,\kappa))$ is the largest root of the equation:
\begin{equation}\label{e2.6}
f(\lambda ; a,b,\kappa):=\lambda^3 - (n - 3)\lambda^2 + (ab-2n+3)\lambda +ab(\kappa+1)-n+1=0.
\end{equation}
%Let $f(a,b,\kappa;\lambda)=\lambda^3 - (n - 3)\lambda^2 + (ab-2n+3)\lambda +ab(\kappa+1)-n+1$.
Then, we have
$$f(\lambda ; a,b,\kappa)-f(\lambda ; \delta-\kappa+1,n-\delta-1,\kappa)=(\lambda+\kappa+1)(ab-(\delta-\kappa+1)(n-\delta-1))\geq 0$$
for any $\lambda>0$ and $b\geq a\geq \delta-\kappa+1$.
Therefore, $\rho(G(a,b,\kappa))\leq \rho(G(\delta-\kappa+1,n-\delta-1,\kappa))$ for any $b\geq a\geq \delta-\kappa+1$,
which means
$$\rho(K_{\kappa}\vee (K_{|X_0|}\cup K_{|Y|}))\leq \rho(K_{\kappa}\vee (K_{\delta-\kappa+1}\cup K_{n-\delta-1})).$$
Since $K_{\kappa}\vee (K_{\delta-\kappa+1}\cup K_{n-\delta-1})$ is a subgraph of $K_{k-1}\vee (K_{\delta-k+2}\cup K_{n-\delta-1})$
for any $\kappa\leq k-1$,
$$\rho(K_{\kappa}\vee (K_{\delta-\kappa+1}\cup K_{n-\delta-1}))\leq \rho(K_{k-1}\vee (K_{\delta-k+2}\cup K_{n-\delta-1})).$$
Hence, from the discussion above we have
$$
\rho(G)\leq \rho(K_{\kappa}\vee (K_{|X_0|}\cup K_{|Y|}))\leq \rho(K_{\kappa}\vee (K_{\delta-\kappa+1}\cup K_{n-\delta-1}))\leq \rho(K_{k-1}\vee (K_{\delta-k+2}\cup K_{n-\delta-1})).
$$
By (\ref{e2.5}), the inequalities above must be equalities. Thus $|X_0|=\delta-k+2$, $\kappa=k-1$, $|Y|=n-\delta-1$
and $G = K_{k-1}\vee (K_{\delta-k+2}\cup K_{n-\delta-1})$.
The proof is completed.
\end{pf}

\begin{rem}{\em
In the Corollary 3.5 of reference \cite{fzl17}, the authors showed that if $G$ is a connected graph of 
minimum degree $\delta(G)\geq \delta\geq k\geq 3$ and order $n\geq (\delta-k+2)(k^2-2k+4)+3$, 
and $\rho(G) \geq \rho(K_{k-1}\vee (K_{\delta-k+2}\cup K_{n-\delta-1}))$,
then $G$ is $k$-connected unless $G = \rho(K_{k-1}\vee (K_{\delta-k+2}\cup K_{n-\delta-1}))$.
Apparently, Theorem \ref{thm2.2} improves the Corollary 3.5 of reference \cite{fzl17} without restriction
on the order of graph. }
\end{rem}

\begin{thm}\label{thm2.3}
Let $G$ be a connected graph of order $n$, minimum degree $\delta\geq k\geq 3$. If
$G$ is a subgraph of $K_{k-1}\vee (K_{\delta-k+2}\cup K_{n-\delta-1})$ and $n\geq \frac{1}{2}(\delta-k+2)(k^2-2k+7)$, then
\begin{equation*}%\label{ }
\rho(G) < n-\delta+k-3,
\end{equation*}
unless $G=K_{k-1}\vee (K_{\delta-k+2}\cup K_{n-\delta-1})$.
\end{thm}

\begin{pf}
In short, denote $H=K_{k-1}\vee (K_{\delta-k+2}\cup K_{n-\delta-1})$.
Let $\textbf{x} = (x_1, x_2,\dots , x_n)^T$ be the unique positive unit eigenvector
corresponding to $\rho(G)$. Recalling that Rayleigh's principle implies that
$$\rho(G)=\textbf{x}^TA(G)\textbf{x}=2\sum_{ij\in E(G)}x_ix_j.$$

Assume that $G$ is a proper subgraph of $H$. Clearly, we could
assume that $G$ is obtained by omitting just one edge $uv$ of $H$.
Let $X, Y, Z$ be the set of vertices of $H$ of degree $\delta, n-1, n-\delta+k-3$, respectively,
where $|X|=\delta-k+2$, $|Y|=k-1$, $|Z|=n-\delta-1$.
Since $\delta(G)=\delta$, $G$ must contain all the edges between $X$ and $Y$.
Therefore, $\{u,v\} \subset Y\cup Z$, with three possible cases: (a) $\{u,v\}\subset Y$ ; (b) $u \in Y, v \in Z$;
(c) $\{u,v\}\subset Z$. % Note that if $k=2$, then case (a) can not happen.
We shall show that case (c) yields a graph of no smaller spectral radius
than case (b), and that case (b) yields a graph of no smaller spectral radius than
case (a).

Firstly, suppose that case (a) occurs, that is, $\{u,v\}\subset Y$. Choose a vertex $w\in Z$.
If $x_u\geq x_w$, then by removing the edge $vw$ and adding the edge $uv$
we obtain a new graph $G_1$ which is covered by case (b). By the Rayleigh principle,
$$\rho(G_1)-\rho(G)\geq \textbf{x}^TA(G_1)\textbf{x}-\textbf{x}^TA(G)\textbf{x}=2x_v(x_u-x_w)\geq 0.$$
If $x_w> x_u$, then by removing all the edges
between $X$ and $\{u\}$ and adding all the edges between $X$ and $\{w\}$
we obtain a new graph $G_1'$ which is also covered by case (b). By the Rayleigh principle,
$$\rho(G_1')-\rho(G)\geq \textbf{x}^TA(G_1')\textbf{x}-\textbf{x}^TA(G)\textbf{x}=2(x_w-x_u)\sum\limits_{i\in X}x_i> 0.$$

Secondly, suppose that case (b) occurs, that is, $u\in Y, v\in Z$. Choose a vertex $w\in Z$ and $w\neq v$.
If $x_u\geq x_w$, then by removing the edge $vw$ and adding the edge $uv$
we obtain a new graph $G_2$ which is covered by case (c). By the Rayleigh principle,
$$\rho(G_2)-\rho(G)\geq \textbf{x}^TA(G_2)\textbf{x}-\textbf{x}^TA(G)\textbf{x}=2x_v(x_u-x_w)\geq 0,$$
If $x_w> x_u$, then by removing all the edges
between $X$ and $\{u\}$ and adding all the edges between $X$ and $\{w\}$
we obtain a new graph $G_2'$ which is also covered by case (c). By the Rayleigh principle,
$$\rho(G_2')-\rho(G)\geq \textbf{x}^TA(G_2')\textbf{x}-\textbf{x}^TA(G)\textbf{x}=2(x_w-x_u)\sum\limits_{i\in X}x_i> 0.$$

Therefore, we may assume that $\{u,v\} \in Z$. By symmetry,
let $x:= x_i$ for any $i\in X$; $y:= x_j$ for any
$j\in Y$; $z := x_\ell $ for any $\ell \in Z\setminus\{u, v\}$; and $t := x_u = x_v$.
According to $\lambda x_i=\sum\limits_{ij\in E(G)}x_j$ and the uniqueness of $\textbf{x}$, we have that
$\rho$ is the largest root of following equations:
\begin{align*}
\lambda x &=(\delta-k+1)x+(k-1)y,\\
\lambda y &=(\delta-k+2)x+(k-2)y+(n-\delta-3)z+2t,\\
\lambda z &=(k-1)y+(n-\delta-4)z+2t,\\
\lambda t &=(k-1)y+(n-\delta-3)z.
\end{align*}
Thus $\rho(G)$ is the largest root of the equation:
\begin{align*}
f(\lambda):= &\lambda^4 - (n - 5)\lambda^3 + ((n-\delta-1)(\delta-k-2) - 4\delta+7)\lambda^2 \\
 &+ [(\delta k+2\delta+2)( n - \delta+k-3) - (k^2+3)(n - 1)+6]\lambda \\
 &+ 2((\delta-k+1)(k(n-\delta-2)-1)+(k-1)(n-\delta-3))=0.
\end{align*}
%%%  hh=x^4 - (n - 5)*x^3 + ((n-d-1)*(d-k-2) - 4*d+7)*x^2 + ((d*k+2*d+2)*( n - d+k-3) - (k^2+3)*(n - 1)+6)*x + 2*((d-k+1)*(k*(n-d-2)-1)+(k-1)*(n-d-3))
By some basic calculations, we have
$$f(n-\delta+k-3)=2n^2-(\delta-k+2)(k^2-2k+7)n+(\delta-k+2)((\delta-k+1)(k^2-2k+5)-2(k-3)).$$
Set $g(x)=2x^2 -(\delta-k+2)(k^2-2k+7)x+(\delta-k+2)((\delta-k+1)(k^2-2k+5)-2(k-3))$.
It is easy to see that the function $g(x)$ is strictly increasing when
$x>\frac14(\delta-k+2)(k^2-2k+7)$. Since $n\geq \frac12(\delta-k+2)(k^2-2k+7) > \frac14(\delta-k+2)(k^2-2k+7)$,
we can get
\begin{align*}
f(n-\delta+k-3)=g(n)\geq & g(\tfrac12(\delta-k+2)(k^2-2k+7))\\
= & (\delta-k+2)((\delta-k+1)(k^2-2k+5)-2(k-3))\\
\geq & (\delta-k+2)(k^2-4k+11)>0.
\end{align*}
By $3\leq k\leq \delta\leq n-3$ and $n\geq \frac{1}{2}(\delta-k+2)(k^2-2k+7)$, we have $n\geq 2(\delta-k+3)$ and so
\begin{align*}
f(n-\delta+k-4)= & - n^3 +4(\delta - k + 3)n^2-((k^2-8k+5\delta+23)(\delta-k+2)-5)n \\
                 &  + ((k^2 - 5k + 2\delta + 15)(\delta - k + 2)-(2\delta - 5))(\delta-k+2)+2  \\
               = &-n(n-2(\delta-k+3))^2-((k^2-4k+\delta+7)(\delta-k+2)+1)(n-2(\delta-k+3)) \\
                 &- (\delta - k + 2)((k^2 - 3k + 3)(\delta - k + 3)+ k(k - 1)) \\
               \leq & -2(3\cdot 3+k(k-1))<0,
\end{align*}
$f(0) =2((\delta-k+1)(k(n-\delta-2)-1)+(k-1)(n-\delta-3))\geq 2(k-1)>0$,  % $f(-\sqrt{3}) =2\sqrt{3}-4<0$ when $\delta=k=2$,
$f(-2) =-2(k - 2)(\delta - k + 2)+2<0$, and $f(-\infty) >0$.
%Therefore, there is at least one root of $f(x) = 0$ between $n-3$ and $n-4$. That is, there is
%at least one eigenvalue of $A(G)$ between $n-4$ and $n-3$.
Therefore, it is easy to find that the largest root of $f(x) = 0$ is in the interval $(n-\delta+k-4,n-\delta+k-3)$,
and it follows that $\rho(G) <n-\delta+k-3$.
\end{pf}

\begin{lem}{\rm (Hong et al.~\cite{h01})}\label{lem2.4}
Let $G$ be a connected graph with $n$ vertices and $m$ edges. Let
$\delta=\delta(G)$ be the minimum degree of $G$ and $\rho(G)$ be the spectral
radius of the adjacency matrix of $G$. Then
$$
\rho(G)\leq \frac{\delta-1}{2}+\sqrt{2m-\delta n+\frac{(\delta+1)^2}{4}}.
$$
Equality holds if and only if $G$ is either a regular graph or a bidegreed graph
in which each vertex is of degree either $\delta$ or $n-1$.
\end{lem}

\begin{thm}\label{thm2.5}
Let $G$ be a connected graph of order $n$, minimum degree $\delta\geq k\geq 3$ and vertex-connectivity $\kappa$. If
$n\geq \frac12(\delta-k+2)(k^2-2k+7)$ and
\begin{equation*}%\label{}
\rho(G) \geq n-\delta+k-3,
\end{equation*}
then $\kappa\geq k$, unless $G=K_{k-1}\vee (K_{\delta-k+2}\cup K_{n-\delta-1})$.
\end{thm}

\begin{pf}
On the contrary, suppose that $\kappa<k$.
Since $G$ is connected and $\rho(G) \geq n-\delta+k-3$, by Lemma~\ref{lem2.4}, we have
\begin{equation*}%\label{}
n-\delta+k-3  \leq \rho(G) \leq \frac{\delta-1}{2}+\sqrt{2|E(G)|-\delta n+\frac{(\delta+1)^2}{4}},
\end{equation*}
which yields
$$|E(G)|\geq \frac{1}{2}n(n-1)-\frac{1}{2}(\delta-k+2)(2n-2\delta+k-3).$$
Since $n\geq \frac12(\delta-k+2)(k^2-2k+7)$, we obtain $n\geq \frac12(k+1)(\delta-k+2)+(\delta+2)$.
By Theorem~\ref{thm2.1} (b), $G$ is a subgraph of $K_{k-1}\vee (K_{\delta-k+2}\cup K_{n-\delta-1})$.
Since $\rho(G) \geq n-\delta+k-3$, by Theorem~\ref{thm2.3}, $G=K_{k-1}\vee (K_{\delta-k+2}\cup K_{n-\delta-1})$.
The proof is completed.
\end{pf}

%\vspace{10pt}

\begin{rem}{\em
In the Theorem 3.4 of reference \cite{fzl17}, the authors proved that if $G$ is a connected graph of
minimum degree $\delta(G)\geq \delta\geq k\geq 3$ and order $n\geq (\delta-k+2)(k^2-2k+4)+3$,
and $\rho(G) \geq n-\delta+k-3$,
then $G$ is $k$-connected unless $G = \rho(K_{k-1}\vee (K_{\delta-k+2}\cup K_{n-\delta-1}))$.
Obviously, Theorem \ref{thm2.5} improves the Theorem 3.4 of reference \cite{fzl17} for the restriction
on the order of graph. }
\end{rem}

Another sufficient condition for graphs to be $k$-connected can be obtained by using the spectral radius of the complement of a graph.

\begin{thm}\label{thm2.6}
Let $G$ be a connected graph of order $n\geq 5$, minimum degree $\delta\geq k\geq 2$ and vertex-connectivity $\kappa$. If
\begin{equation}\label{e2.7}
\rho(\overline{G}) \leq \sqrt{(\delta-k+2)(n-\delta-1)},
\end{equation}
then $\kappa\geq k$, unless $G=K_{k-1}\vee (K_{\delta-k+2}\cup K_{n-\delta-1})$.
\end{thm}

\begin{pf}
%For short, let $H=K_{\delta-1}\vee (K_2\cup K_{n-\delta-1})$.
Assume that (\ref{e2.7}) holds but $1\leq \kappa\leq k-1$.
Let $S$ be an arbitrary minimum vertex-cut of $G$, and
let $X_0, X_1,\dots$, $X_{p-1}$ ($p\geq 2$) denote the vertex sets of the components of $G-S$, where
$|X_0|\leq |X_1|\leq \cdots\leq |X_{p-1}|$. Each vertex in
$X_i$ is adjacent to at most $|X_i|-1$ vertices of $X_i$ and $\kappa=|S|$ vertices of $S$. Thus
\begin{equation*}
\delta|X_i| \leq \sum_{x\in X_i}d(x) \leq |X_i|(|X_i|-1+\kappa) ,
\end{equation*}
and so
$|X_i|\geq \delta-\kappa+1$ for each $i=0,1,\dots,p-1$.
Let $Y = \bigcup_{i=1}^{p-1} X_i$. Then $\delta-\kappa+1 \leq |X_0|\leq |Y| \leq n-\delta-1$ and $|X_0|+|Y| = n-\kappa$.
Since there are no edges between $X_0$ and $Y$ in $G$, $K_{|X_0|,|Y|}$ is a subgraph of $\overline{G}$.
Thus
\begin{align*}
\rho(\overline{G}) & \geq  \rho(K_{|X_0|,|Y|})=\sqrt{|X_0|\cdot|Y|}= \sqrt{|X_0|(n-\kappa-|X_0|)}\\
& \geq  \sqrt{(\delta-\kappa+1)(n-\delta-1)}\geq \sqrt{(\delta-k+2)(n-\delta-1)}.
\end{align*}
By (\ref{e2.7}), the inequalities above must be equalities. Thus $|X_0|=\delta-k+2$, $\kappa=k-1$
and $\overline{G}=K_{\delta-k+2,n-\delta-1}$, and so $G=K_{k-1}\vee (K_{\delta-k+2}\cup K_{n-\delta-1})$.
\end{pf}

%Since $\kappa\leq \delta-1$ and $|X_0|\geq 2$, $K_{2,n-\delta-1}\subseteq \overline{G}$. Thus, we have
%$$
%\rho(\overline{G})\geq \rho(K_{2,n-\delta-1})=\sqrt{2(n-\delta-1)}.
%$$

\section{Maximally connected graphs}

If $\kappa(G)=\delta(G)$, then $G$ is maximally connected. Therefore, by setting $k=\delta$ in Theorem \ref{thm2.1},
we obtain the following theorem.

\begin{thm}\label{thm3.1}
Let $G$ be a connected graph of order $n\geq 5$, size $m$, minimum degree $\delta\geq 2$ and vertex-connectivity $\kappa$.

{\rm (a)}\
If
%\begin{equation*}%\label{e3.1}
$m \geq {n-2 \choose 2}+2\delta-1$,
%\end{equation*}
then $\kappa=\delta$, unless $G=K_{\delta-1}\vee (K_2\cup K_{n-\delta-1})$.

{\rm (b)}\
If $n\geq  2\delta+3$ and
%\begin{equation*}%\label{e3.2}
$m \geq {n-2 \choose 2}+\delta$,
%\end{equation*}
then $\kappa=\delta$, unless $G$ is a subgraph of $K_{\delta-1}\vee (K_2\cup K_{n-\delta-1})$.
\end{thm}

\begin{thm}\label{thm3.2}
Let $G$ be a connected graph of order $n\geq 5$, minimum degree $\delta\geq 2$ and vertex-connectivity $\kappa$. If
\begin{equation}\label{e3.3}
\rho(G) \geq \frac{\delta-1}{2}+\sqrt{(n-\delta-1)(n-4) + \frac{(\delta+1)^2}{4}},
\end{equation}
then $\kappa=\delta$, unless $G=K_{n-4}\vee (K_2\cup K_2)$.
\end{thm}

\begin{pf}
On the contrary, suppose that $\kappa<\delta$.
Since $G$ is connected, by (\ref{e3.3}) and Lemma~\ref{lem2.4}, we have
\begin{equation}\label{e3.4}
\frac{\delta-1}{2}+\sqrt{(n-\delta-1)(n-4) + \frac{(\delta+1)^2}{4}}  \leq \rho(G) \leq \frac{\delta-1}{2}+\sqrt{2|E(G)|-\delta n+\frac{(\delta+1)^2}{4}},
\end{equation}
which yields
$$|E(G)|\geq {n-2 \choose 2}+2\delta-1.$$
By Theorem~\ref{thm3.1} (a), $G=K_{\delta-1}\vee (K_2\cup K_{n-\delta-1})$.
To complete the proof, we only need to show $\delta=n-3$.

Since $|E(G)|={n-2 \choose 2}+2\delta-1$, the equalities hold in (\ref{e3.4}).
Thus, by Lemma~\ref{lem2.4}, $G$ is a bidegreed graph in which each vertex is of degree $\delta$ or $n-1$.
However, the vertices of $G$ have degrees from the set $\{\delta,n-3,n-1\}$.
Therefore, $\delta=n-3$ and the result follows.
\end{pf}

\vspace{9pt}
By setting $k=\delta$ in Theorem \ref{thm2.2}, we obtain the following result.

\begin{thm}\label{thm3.3}
Let $G$ be a connected graph of order $n\geq 5$, minimum degree $\delta\geq 2$ and vertex-connectivity $\kappa$. If
\begin{equation*}\label{e3.5}
\rho(G) \geq \rho(K_{\delta-1}\vee (K_2\cup K_{n-\delta-1})) ,
\end{equation*}
then $\kappa=\delta$, unless $G=K_{\delta-1}\vee (K_2\cup K_{n-\delta-1})$, where $\rho(K_{\delta-1}\vee (K_2\cup K_{n-\delta-1}))$
is the largest root of the equation
$$\lambda^3 - (n - 3)\lambda^2 - (2\delta-1)\lambda + 2\delta(n-\delta-1)-n+1=0.$$
\end{thm}

\begin{thm}\label{thm3.4}
Let $G$ be a connected graph of order $n$, minimum degree $\delta\geq 2$ and vertex-connectivity $\kappa$. If
$n\geq \delta^2-2\delta+7$ and
\begin{equation*}%\label{e1.6}
\rho(G) \geq n-3,
\end{equation*}
then $\kappa=\delta$, unless $G=K_{\delta-1}\vee (K_2\cup K_{n-\delta-1})$.
\end{thm}

\begin{pf}
Set $k=\delta$ in the proofs of Theorem~\ref{thm2.3} and Theorem~\ref{thm2.5}.
If $\delta\geq 3$, then the result follows by Theorem~\ref{thm2.5}.
If $\delta=2$, then case (a) can not occur in the proof of Theorem~\ref{thm2.3}.
In Theorem~\ref{thm2.3}, by noting that $f(n-3)>0$, $f(n-4)<0$, $f(0)>0$, $f(-\sqrt{3})=2\sqrt{3}-4<0$ and $f(-\infty) >0$,
it follows that $\rho(G)<n-3$ and so Theorem~\ref{thm2.3} holds for $\delta=k=2$.
Hence, Theorem~\ref{thm2.5} also holds for $\delta=k=2$ and the result follows.
\end{pf}

\begin{rem}{\rm
In the proof of Theorem~\ref{thm2.3}, if we take $k=\delta\geq 2$ and $n = \delta^2-2\delta+6$, then
$f(n-3)=g(n) = g(\delta^2-2\delta+6)=10-4\delta < 0$
when $\delta \geq 3$. Notice that $f(+\infty) = +\infty$. So, the largest root of $f(x) = 0$ is greater than $n-3$
if $\delta \geq 3$, and it follows that $\rho(G) > n -3$. That is to say, the requirement $n\geq\delta^2-2\delta+7$
in Theorem \ref{thm3.4} is best possible when $\delta \geq 3$. }
\end{rem}

By setting $k=\delta$ in Theorem \ref{thm2.6}, we have the following result.

\begin{thm}\label{thm3.6}
Let $G$ be a connected graph of order $n\geq 5$, minimum degree $\delta\geq 2$ and vertex-connectivity $\kappa$. If
\begin{equation*}%\label{e3.6}
\rho(\overline{G}) \leq \sqrt{2(n-\delta-1)},
\end{equation*}
then $\kappa=\delta$, unless $G=K_{\delta-1}\vee (K_2\cup K_{n-\delta-1})$.
\end{thm}

\section{Super-connected graphs}

For any connected graph $G$ of order $n$,
if $2\leq n\leq 4$, then $G$ is super-$\kappa$. Therefore,
$n\geq 5$ is considered in this section.

\begin{thm}\label{thm4.1}
Let $G$ be a connected graph of order $n\geq 5$, size $m$, minimum degree $\delta$ and vertex-connectivity $\kappa$. If
\begin{equation}\label{e4.1}
m \geq {n-2 \choose 2}+2\delta,
\end{equation}
then $G$ is super-$\kappa$, unless $G=(K_{\delta}\vee (K_2\cup K_{n-\delta-2}))-e$, where $e=xy$ is an edge of $K_{\delta}\vee (K_2\cup K_{n-\delta-2})$
with $d(x)=\delta+1$ and $d(y)=n-1$.
\end{thm}

\begin{pf}
Since $m \geq {n-2 \choose 2}+2\delta$, by Theorem~\ref{thm3.1} (a), $\kappa=\delta$.
%If $\delta\geq n-3$, then $\kappa\geq n-3$ and so $G$ is super-$\kappa$.
%Thus assume in the following that $\delta\leq n-4$.
On the contrary, suppose that $G$ is not super-$\kappa$.
Let $S$ be an arbitrary minimum vertex-cut with $\delta$ vertices,
and let $X_0, X_1,\dots, X_{p-1}$ ($p\geq 2$) denote the vertex sets of the components of $G-S$, where
$2\leq |X_0|\leq |X_1|\leq \cdots\leq |X_{p-1}|$. Denote $Y = \bigcup_{i=1}^{p-1} X_i$.
Since $G-S$ is disconnected, there are no edges between $X_0$ and $Y$ in $G$, and
\begin{equation}\label{e4.2}
m \leq \frac12 n(n-1) - |X_0|\cdot|Y|.
\end{equation}
Thus, by $|X_0|+|Y|=n-\delta$ and $2 \leq |X_0|\leq |Y| \leq n-\delta-2$, we have
\begin{align*}
m \leq \frac12 n(n-1) - |X_0|\cdot|Y|
  \leq \frac12 n(n-1) - 2(n-\delta-2)
  ={n-2 \choose 2}+2\delta+1.
\end{align*}

If $m= {n-2 \choose 2}+2\delta+1$, then all the inequalities in the proof above must be equalities.
We have $p=2$, $|X_0|=2$, $|Y|=n-\delta-2$, $d_G(s)=n-1$ for each $s\in S$,
$d_{G[X_0]}(x)=1$ for each $x\in X_0$ and $d_{G[Y]}(y)=n-\delta-3$ for each $y\in Y$. That is
$G[X_0]=K_2$, $G[S]=K_{\delta}$, $G[Y]=K_{n-\delta-2}$ and $G\cong K_{\delta}\vee (K_2\cup K_{n-\delta-2})$.
However, $\delta(G)=\delta+1>\delta$, which is a contradiction. Therefore, $m\leq {n-2 \choose 2}+2\delta$.
By (\ref{e4.1}), $m = {n-2 \choose 2}+2\delta$.

Next, we show that $G$ is a proper subgraph of $K_{\delta}\vee (K_2\cup K_{n-\delta-2})$.
It suffices to prove that $|X_0|=2$ and $p=2$.

If $|X_0|\geq 3$, then $n\geq \delta+6$. Combining (\ref{e4.1}) with (\ref{e4.2}), we obtain
\begin{align*}
{n-2 \choose 2}+2\delta \leq m &\leq \frac{1}{2}(n^2 -n)- |X_0|\cdot|Y|\\
                                 &\leq \frac{1}{2}(n^2 -n)-3(n-\delta-3) \\
                                 &\leq {n-2 \choose 2}+2\delta.
\end{align*}
All the inequalities above must be equalities, and so $|X_0|=|Y|= 3$, $n=\delta+6$, $G\cong K_{\delta}\vee (K_3\cup K_3)$.
However, $\delta(G)=\delta+2>\delta$, which is a contradiction. Therefore, $|X_0|=2$.

If $p\geq 3$, then $n\geq \delta+6$. Let $Y_1=\bigcup_{i=2}^{p-1} X_i$.
Then $|X_1| + |Y_1|=n-\delta-2$ and $2=|X_0|\leq |X_1|\leq |Y_1|\leq n-\delta-4$.
Since $G-S$ is disconnected, there are no edges among $X_0$, $X_1$ and $Y_1$ in $G$
(i.e. $[X_0,X_1]=\emptyset$, $[X_0,Y_1]=\emptyset$, $[X_1,Y_1]=\emptyset$), and
\begin{align*}
m  & \leq \frac{1}{2}(n^2 -n)- |X_0|\cdot(|X_1|+|Y_1|)-|X_1||Y_1|\\
   & \leq \frac{1}{2}(n^2 -n)- 2(n-\delta-2)-2\cdot 2\\
   & = {n-2 \choose 2}+2\delta-3  < {n-2 \choose 2}+2\delta,
\end{align*}
which is a contradiction. Therefore, $p=2$.

Let $H=K_{\delta}\vee (K_2\cup K_{n-\delta-2})$. Then $G\subset H$  and
$|E(G)|=|E(H)|-1$. Therefore, $G=H-e$.
Since $\delta(H)=\delta+1$ and $\delta(G)=\delta$,
$e=xy$ is an edge of $H$ with $d(x)=\delta+1$ and $d(y)=n-1$.
\end{pf}

\begin{thm}\label{thm4.2}
Let $G$ be a connected graph of order $n\geq 5$, minimum degree $\delta$. If
\begin{equation}\label{e4.3}
\rho(G) \geq \frac{\delta-1}{2}+\sqrt{2+(n-\delta-1)(n-4) + \frac{(\delta+1)^2}{4}},
\end{equation}
then $G$ is super-$\kappa$.
\end{thm}

\begin{pf}
On the contrary, suppose that $G$ is not super-$\kappa$.
Since $G$ is connected, by (\ref{e4.3}) and Lemma~\ref{lem2.4}, we have
%$$\sqrt{(n-1)(n-5) + 4\delta+2} \leq \rho(G) \leq \sqrt{2|E(G)|-n+1},$$
\begin{equation}\label{e3.6}
\sqrt{2+(n-\delta-1)(n-4) + \frac{(\delta+1)^2}{4}} \leq \rho(G)-\frac{\delta-1}{2} \leq \sqrt{2|E(G)|-\delta n+\frac{(\delta+1)^2}{4}},
\end{equation}
which yields
$$|E(G)|\geq {n-2 \choose 2}+2\delta.$$
By Theorem~\ref{thm4.1},
$G=K_{\delta}\vee (K_2\cup K_{n-\delta-2})-e$, where $e=xy$ is an edge of $K_{\delta}\vee (K_2\cup K_{n-\delta-2})$
with $d(x)=\delta+1$ and $d(y)=n-1$.

Since $|E(G)|={n-2 \choose 2}+2\delta$, the equalities hold in (\ref{e3.6}).
Thus, by Lemma~\ref{lem2.4}, $G$ is a bidegreed graph in which each vertex is of degree $\delta$ or $n-1$.
However, the vertices of $G$ have degree from the set
$\{\delta,\delta+1,n-3,n-1\}$. Thus $G$ cannot be a bidegreed graph, which yields a contradiction.
Hence, $G$ is super-$\kappa$.
\end{pf}

\begin{thm}\label{thm4.3}
Let $G$ be a connected graph of order $n\geq 5$, minimum degree $\delta$. If
\begin{equation}\label{e4.4}
\rho(G) \geq \rho(K_\delta\vee (K_2\cup K_{n-\delta-2})),
\end{equation}
then $G$ is super-$\kappa$, where $\rho(K_{\delta}\vee (K_2\cup K_{n-\delta-2}))$
is the largest root of the equation
$$\lambda^3 - (n - 3)\lambda^2 - (2\delta+1)\lambda + 2\delta(n-\delta-2)-n+1=0.$$
\end{thm}

\begin{pf}
On the contrary, suppose that $G$ is not super-$\kappa$.
Let $S$ be an arbitrary minimum vertex-cut with $\kappa$ $(\leq\delta)$ vertices,
and let $X_0, X_1,\dots, X_{p-1}$ ($p\geq 2$) denote the vertex sets of the components of $G-S$, where
$2\leq |X_0|\leq |X_1|\leq \cdots\leq |X_{p-1}|$. Denote $Y = \bigcup_{i=1}^{p-1} X_i$.
Then $2\leq |X_0|\leq |Y|\leq n-\kappa-2$ and $|X_0|+|Y|=n-\kappa$.
Since there are no edges between $X_0$ and $Y$ in $G$,
$G$ is a subgraph of $K_{\kappa}\vee (K_{|X_0|}\cup K_{|Y|})$ and
$\rho(G)\leq \rho(K_{\kappa}\vee (K_{|X_0|}\cup K_{|Y|}))$.

According to (\ref{e2.6}) in the proof of Theorem \ref{thm2.2},
$\rho(K_{\kappa}\vee (K_{|X_0|}\cup K_{|Y|}))$ is the largest root of the equation:
$$f(\lambda ; |X_0|,|Y|,\kappa):=\lambda^3 - (n - 3)\lambda^2 + (|X_0||Y|-2n+3)\lambda +|X_0||Y|(\kappa+1)-n+1=0.$$
Then, we have
$$f(\lambda ; |X_0|,|Y|,\kappa)-f(\lambda ; 2,n-\kappa-2,\kappa)=(\lambda+\kappa+1)(|X_0||Y|-2(n-\kappa-2))\geq 0$$
for any $\lambda>0$ and $|Y|\geq |X_0|\geq 2$.
Therefore, $\rho(K_{\kappa}\vee (K_{|X_0|}\cup K_{|Y|}))\leq \rho(K_{\kappa}\vee (K_{2}\cup K_{n-\kappa-2}))$.

Since $K_{\kappa}\vee (K_{2}\cup K_{n-\kappa-2})$ is a subgraph of $K_{\delta}\vee (K_{2}\cup K_{n-\delta-2})$
for any $\kappa\leq \delta$, we get
$$\rho(K_{\kappa}\vee (K_{2}\cup K_{n-\kappa-2}))\leq \rho(K_{\delta}\vee (K_{2}\cup K_{n-\delta-2})).$$

Hence, from the discussion above we obtain
$$
\rho(G)\leq \rho(K_{\kappa}\vee (K_{|X_0|}\cup K_{|Y|}))\leq \rho(K_{\kappa}\vee (K_{2}\cup K_{n-\kappa-2}))\leq \rho(K_{\delta}\vee (K_{2}\cup K_{n-\delta-2})).
$$
By (\ref{e4.4}), the above inequalities must be equalities. Thus $|X_0|= 2$, $\kappa=\delta$, $|Y|=n-\delta-2$
and $G = K_{\delta}\vee (K_{2}\cup K_{n-\delta-2})$. However, $\delta(G)=\delta+1>\delta$, which is a contradiction.
The result follows.
\end{pf}

\begin{thm}\label{thm4.4}
Let $G$ be a connected graph of order $n\geq 5$, minimum degree $\delta$ and vertex-connectivity $\kappa$. If
\begin{equation}\label{e4.5}
\rho(\overline{G}) \leq \sqrt{2(n-\delta-2)},
\end{equation}
then $G$ is super-$\kappa$.
\end{thm}

\begin{pf}
%If $\delta\geq n-3$, then $\kappa\geq n-3$ and so $G$ is super-$\kappa$.
%Thus assume in the following that $\delta\leq n-4$.
%Assume that (\ref{e4.5}) holds but $G$ is not super-$\kappa$.
%If $\kappa\leq \delta-1$, then, by the proof of Theorem \ref{thm2.6} and $\delta\leq n-4$,
%$$
%\rho(\overline{G})\geq \sqrt{2(n-\delta-1)}>\sqrt{n-\delta-\frac{3}{2}+\sqrt{(n-\delta-2)^2+\frac{1}{4}}},
%$$
%a contradiction. Next, assume that $\kappa = \delta$ in the following.
%Let $S$ be an arbitrary minimum vertex-cut of $G$ with $\kappa=\delta$ vertices, and
%let $X_0, X_1,\dots$, $X_{p-1}$ ($p\geq 2$) denote the vertex sets of the components of $G-S$, where
%$2\leq |X_0|\leq |X_1|\leq \cdots\leq |X_{p-1}|$.
Assume that (\ref{e4.5}) holds but $G$ is not super-$\kappa$.
Let $S$ be an arbitrary minimum vertex-cut of $G$ with $\kappa$ $(\leq \delta)$ vertices, and
let $X_0, X_1,\dots$, $X_{p-1}$ ($p\geq 2$) denote the vertex sets of the components of $G-S$, where
$2\leq |X_0|\leq |X_1|\leq \cdots\leq |X_{p-1}|$.
Denote $Y = \bigcup_{i=1}^{p-1} X_i$. Then $2 \leq |X_0|\leq |Y| \leq n-\kappa-2$ and $|X_0|+|Y| = n-\kappa$.
Since there are no edges between $X_0$ and $Y$ in $G$, $K_{|X_0|,|Y|}$ is a subgraph of $\overline{G}$.
Thus
$$
\rho(\overline{G})\geq \rho(K_{|X_0|,|Y|})=\sqrt{|X_0|\cdot|Y|}= \sqrt{|X_0|(n-\kappa-|X_0|)}\geq \sqrt{2(n-\kappa-2)}\geq \sqrt{2(n-\delta-2)}.
$$
By (\ref{e4.5}), the inequalities above must be equalities. Thus $|X_0|=2$, $\kappa=\delta$
and $\overline{G}=K_{2,n-\delta-2}$, and so $G=K_{\delta}\vee (K_2\cup K_{n-\delta-2})$.
However, $\delta(G)=\delta(K_{\delta}\vee (K_2\cup K_{n-\delta-2}))=\delta+1>\delta$, a contradiction.
This completes the proof.
\end{pf}

\section{Triangle-free graphs}

Let us extend an interesting result by applying the famous theorem of Mantel \cite{mantel} and Tur\'{a}n \cite{turan}.

\begin{thm}\label{mantel}
For any triangle-free graph $G$ of order $n$, we have
$|E(G)|\leq \lfloor\frac{1}{4}n^2\rfloor$, with equality if and only if $G=K_{\lfloor n/2\rfloor, \lceil n/2\rceil}$.
\end{thm}

\begin{thm}\label{thm5.2}
Let $G$ be a connected triangle-free graph of order $n$, size $m$, minimum degree $\delta\geq k\geq  2$ and vertex-connectivity $\kappa$.
If
\begin{equation}\label{e5.1}
m \geq \delta^2+\left\lfloor\frac14 (n-2\delta+k-1)^2\right\rfloor,
\end{equation}
then $\kappa\geq k$, unless $V(G)=X\cup S\cup Y$ and $S$ is a minimum vertex-cut of $G$ with $G[S]=\overline{K_{k-1}}$,
$G[X\cup S]=K_{\delta, \delta}$ and
$G[Y\cup S] = K_{\lceil(n-2\delta+k-1)/2\rceil,\lfloor(n-2\delta+k-1)/2\rfloor}$.
%$G=K_{\lceil\frac{n+k-1}{2}\rceil,\lfloor\frac{n-k+1}{2}\rfloor}-(K_{\delta,\lfloor\frac{n-k+1}{2}\rfloor-\delta}\cup K_{\delta-k+1,\lceil\frac{n+k-1}{2}\rceil-\delta})$.
%or $G=K_{\lfloor\frac{n+k-1}{2}\rfloor,\lceil\frac{n-k+1}{2}\rceil}-(K_{\delta,\lceil\frac{n-k+1}{2}\rceil-\delta}\cup K_{\delta-k+1,\lfloor\frac{n+k-1}{2}\rfloor-\delta})$.
\end{thm}

\begin{pf}
On the contrary, suppose that $\kappa\leq k-1$. Let $S$ be a minimum vertex-cut of $G$, and let
$X, X_1,\dots, X_{p-1}$ ($p\geq 2$) denote the vertex sets of the components of $G-S$, where
$|X|\leq |X_1|\leq \cdots\leq |X_{p-1}|$. Set $Y=\bigcup_{i=1}^{p-1} X_i$. Then $|X|\leq |Y|$.
Using Theorem \ref{mantel}, we conclude that
\begin{equation}\label{e5.2}
|E(G[X\cup S])|\leq \left\lfloor\frac{(|X|+|S|)^2}{4}\right\rfloor \text{~and~} |E(G[Y\cup S])|\leq \left\lfloor\frac{(|Y|+|S|)^2}{4}\right\rfloor,
\end{equation}
with equalities if and only if
$$G[X\cup S]=K_{\lfloor(|X|+|S|)/2\rfloor, \lceil(|X|+|S|)/2\rceil},
G[Y\cup S]=K_{\lfloor(|Y|+|S|)/2\rfloor, \lceil(|Y|+|S|)/2\rceil}. $$
%Let $(U_1,U_2)$ be a bipartition of $V(G)$, and let $X_i=X\cap U_i$,
%$S_i=S\cap U_i$ and $Y_i=Y\cap U_i$ for $i=1,2$. Thus,
%$X=X_1\cup X_2$, $S=S_1\cup S_2$, $Y=Y_1\cup Y_2$ and $|X|+|Y|=n-\kappa$.
%Note that $|S_1|+|S_2|=\kappa$ and
%$$\delta\leq d_{G}(x_1)\leq |X_2|+|S_2| \text{~for each~} x_1\in X_1, \ \ \delta\leq d_{G}(x_2)\leq |X_1|+|S_1|  \text{~for each~}  x_2\in X_2 ,$$
%$$\delta\leq d_{G}(y_1)\leq |Y_2|+|S_2| \text{~for each~} y_1\in Y_1, \ \ \delta\leq d_{G}(y_2)\leq |Y_1|+|S_1|  \text{~for each~}  y_2\in Y_2 .$$
Note that $|E(G[X\cup S])|\geq \frac{1}{2}\delta(|X|+|S|)$. Thus, by (\ref{e5.2}),
$$
\frac{1}{2}\delta(|X|+|S|)\leq |E(G[X\cup S])|\leq \frac{(|X|+|S|)^2}{4},
$$
and so $|X|\geq 2\delta-|S|=2\delta-\kappa$.
Therefore we arrive at
\begin{equation}\label{e5.3}
2\delta-\kappa \leq |X| \leq |Y|\leq n-2\delta.
\end{equation}
Together with $|X|+|Y|=n-\kappa$ and (\ref{e5.2}), it leads to
\begin{align*}
m & = |E(G[X\cup S])| + |E(G[Y\cup S])| - |E(G[S])| \\
%  & \leq |E(G[X\cup S])| + |E(G[Y\cup S])| \\
  & \leq \left\lfloor\frac{1}{4}(|X|+|S|)^2\right\rfloor + \left\lfloor\frac{1}{4}(|Y|+|S|)^2\right\rfloor - |E(G[S])| \\
  & \leq \left\lfloor\frac{1}{4}(|X|+|S|)^2+\frac{1}{4}(|Y|+|S|)^2\right\rfloor \\
%  & = \left\lfloor\frac{|X|^2+|Y|^2}{4} + \frac{(|X|+|Y|+|S|)|S|}{2} \right\rfloor \\
%  & = \left\lfloor\frac{|X|^2+|Y|^2}{4} + \frac{n\kappa}{2} \right\rfloor \\
  & = \left\lfloor\frac{(|X|+|Y|+|S|)^2+|S|^2}{4} - \frac{|X|\cdot|Y|}{2} \right\rfloor \\
%  & = \left\lfloor\frac{|X|^2+(n-\kappa-|X|)^2}{4} + \frac{n\kappa}{2} \right\rfloor \\
%  & = \left\lfloor\frac{2|X|^2-2(n-\kappa)|X|+(n-\kappa)^2}{4} + \frac{n\kappa}{2} \right\rfloor \\
%  & \leq \left\lfloor\frac{2(2\delta-\kappa)^2-2(n-\kappa)(2\delta-\kappa)+(n-\kappa)^2}{4} + \frac{n\kappa}{2} \right\rfloor \\
%  & \leq \left\lfloor\frac{(2\delta-\kappa)^2+(n-2\delta)^2}{4} + \frac{n\kappa}{2} \right\rfloor \\
   & \leq  \left\lfloor\frac{n^2+\kappa^2}{4} - \frac{(2\delta-\kappa)(n-2\delta)}{2} \right\rfloor \\
%  & = \left\lfloor\delta^2+\frac{1}{4}(n-2\delta+\kappa)^2\right\rfloor \\
  & = \delta^2+ \left\lfloor\frac{1}{4}(n-2\delta+\kappa)^2\right\rfloor \\
  & \leq \delta^2+ \left\lfloor\frac{1}{4}(n-2\delta+k-1)^2\right\rfloor.
\end{align*}
Combining this with (\ref{e5.1}), we have $m =\delta^2+\left\lfloor\frac14 (n-2\delta+k-1)^2\right\rfloor$,
and so $|S|=\kappa=k-1$, $|X|=2\delta-k+1$, $|Y|=n-2\delta$, $|E(G[S])|=0$,
$|E(G[X\cup S])|=\delta^2$ and $|E(G[Y\cup S])| = \left\lfloor\frac14 (n-2\delta+k-1)^2\right\rfloor$.
%and $G[S]=\overline{K_{k-1}}$.
%Since $|U_1|\geq |U_2|$, we can get $S_1=\emptyset$ and $S=S_2$.
%$|X_1|=\delta$, $|X_2|=\delta-k+1$,
%$|Y_1|=\lceil\frac{n+k-1}{2}\rceil-\delta$, $|Y_2|=\lfloor\frac{n-k+1}{2}\rfloor-\delta$.
Therefore, $G[S]=\overline{K_{k-1}}$, $G[X\cup S]=K_{\delta,\delta}$ and $G[Y\cup S] = K_{\lceil(n-2\delta+k-1)/2\rceil,\lfloor(n-2\delta+k-1)/2\rfloor}$.
This completes the proof.
\end{pf}

\vspace{8pt}
By setting $k = \delta$ in Theorem \ref{thm5.2}, we
obtain the following theorem.

\begin{thm}\label{thm5.3}
Let $G$ be a connected triangle-free graph of order $n$, size $m$, minimum degree $\delta\geq  2$ and vertex-connectivity $\kappa$.
If
\begin{equation}\label{e5.4}
m \geq \delta^2+\left\lfloor\frac14 (n-\delta-1)^2\right\rfloor,
\end{equation}
then $\kappa=\delta$, unless $V(G)=X\cup S\cup Y$ and $S$ is a minimum vertex-cut of $G$ with $G[S]=\overline{K_{\delta-1}}$,
$G[X\cup S]=K_{\delta, \delta}$ and $G[Y\cup S] = K_{\lceil(n-\delta-1)/2\rceil,\lfloor(n-\delta-1)/2\rfloor}$.
\end{thm}

\begin{thm}\label{thm5.4}
Let $G$ be a connected triangle-free graph of order $n$, size $m$, minimum degree $\delta\geq  2$ and vertex-connectivity $\kappa$.
If
\begin{equation}\label{e5.5}
m \geq \delta^2+\left\lfloor\frac14 (n-\delta)^2\right\rfloor,
\end{equation}
then $G$ is super-$\kappa$.
%unless $V(G)=X\cup S\cup Y$ and $S$ is a minimum vertex-cut of $G$ with $G[S]=K_{1,\delta-1}$,
%$G[X\cup S]=K_{\delta, \delta}$ and $G[Y\cup S] = K_{\lceil(n-\delta)/2\rceil,\lfloor(n-\delta)/2\rfloor}$.
\end{thm}

\begin{pf}
On the contrary, suppose that $G$ is not super-$\kappa$.
Since $m\geq \delta^2+\left\lfloor\frac14 (n-\delta)^2\right\rfloor$, by Theorem \ref{thm5.3}, $\kappa=\delta$.
Let $S$ be a minimum vertex-cut of $G$ with $\delta$ vertices, and let
$X, X_1,\dots, X_{p-1}$ ($p\geq 2$) denote the vertex sets of the components of $G-S$, where
$2\leq |X|\leq |X_1|\leq \cdots\leq |X_{p-1}|$. Set $Y=\bigcup_{i=1}^{p-1} X_i$. Then $2\leq |X|\leq |Y|$.
%Using Theorem \ref{mantel}, we conclude that
%\begin{equation}\label{e5.6}
%|E(G[X\cup S])|\leq \left\lfloor\frac{(|X|+|S|)^2}{4}\right\rfloor \text{~and~} |E(G[Y\cup S])|\leq \left\lfloor\frac{(|Y|+|S|)^2}{4}\right\rfloor,
%\end{equation}
%with equalities if and only if
%$$G[X\cup S]=K_{\lfloor(|X|+|S|)/2\rfloor, \lceil(|X|+|S|)/2\rceil},
%G[Y\cup S]=K_{\lfloor(|Y|+|S|)/2\rfloor, \lceil(|Y|+|S|)/2\rceil}. $$
%Note that $|E(G[X\cup S])|\geq \frac{1}{2}\delta(|X|+|S|)$. Thus, by (\ref{e5.6}),
%$$
%\frac{1}{2}\delta(|X|+|S|)\leq |E(G[X\cup S])|\leq \frac{(|X|+|S|)^2}{4},
%$$
%and so $|X|\geq 2\delta-|S|=\delta$.
Therefore, with the same proceeding of the proof of Theorem \ref{thm5.2} (from (\ref{e5.2}) to (\ref{e5.3})), we arrive at
\begin{equation}\label{e5.6}
\delta \leq |X| \leq |Y|\leq n-2\delta.
\end{equation}
Together with $|X|+|Y|=n-\delta$ and (\ref{e5.2}), it leads to
\begin{align*}
m & = |E(G[X\cup S])| + |E(G[Y\cup S])| - |E(G[S])| \\
  & \leq \left\lfloor\frac{1}{4}(|X|+|S|)^2\right\rfloor + \left\lfloor\frac{1}{4}(|Y|+|S|)^2\right\rfloor- |E(G[S])| \\
  & \leq \left\lfloor\frac{1}{4}(|X|+|S|)^2+\frac{1}{4}(|Y|+|S|)^2\right\rfloor \\
  & = \left\lfloor\frac{(|X|+|Y|+|S|)^2+|S|^2}{4} - \frac{|X|\cdot|Y|}{2} \right\rfloor \\
  & \leq  \left\lfloor\frac{n^2+\delta^2}{4} - \frac{\delta(n-2\delta)}{2} \right\rfloor \\
  & = \delta^2+ \left\lfloor\frac{1}{4}(n-\delta)^2\right\rfloor .
\end{align*}
Combining this with (\ref{e5.5}), we have $m =\delta^2+\left\lfloor\frac14 (n-\delta)^2\right\rfloor$,
and so $|X|=|S|=\delta$, $|Y|=n-2\delta$, $|E(G[S])|=0$,
$|E(G[X\cup S])|=\delta^2$ and $|E(G[Y\cup S])| = \left\lfloor\frac14 (n-\delta)^2\right\rfloor$.
Therefore, $G[S]=\overline{K_{\delta}}$, $G[X\cup S]=K_{\delta,\delta}$ and $G[Y\cup S] = K_{\lceil(n-\delta)/2\rceil,\lfloor(n-\delta)/2\rfloor}$.
Thus, $G[X]=\overline{K_{\delta}}$, which contradicts to the fact that $G[X]$ is a component of $G$ with at least two vertices.
The result follows.
\end{pf}

\begin{rem}{\rm
The lower bound on $m$ given in Theorem \ref{thm5.4} is sharp. For example,
let $n=3\delta+3$, $V(G)=X\cup S\cup Y$, $G[X]=K_{1,\delta}$, $G[Y]=K_{1,\delta+1}$
and $S$ is a minimum vertex-cut of $G$ with $G[S]=\overline{K_{\delta}}$,
$G[X\cup S]=K_{\delta, \delta+1}$ and $G[Y\cup S] = K_{\delta+1,\delta+1}$.
It is easy to check that
$$|E(G)|=\delta(\delta+1)+(\delta+1)^2=\delta^2+\left\lfloor\frac14 (n-\delta)^2\right\rfloor-1.$$
However, $G-S=K_{1,\delta}\cup K_{1,\delta+1}$, which means $G$ is not super-connected.
}
\end{rem}

\end{document}